\documentclass[12pt]{amsart}
\usepackage{mymacros,xypic,fullpage}

\title{On Images of Mori Dream Spaces}
\author{Shinnosuke Okawa}
\date{\today}


 \DeclareMathOperator{\Nef}{Nef}
 \DeclareMathOperator{\PPic}{\mathbf{Pic}}
 \DeclareMathOperator{\Cl}{Cl}
 \DeclareMathOperator{\Num}{N^1}
 \DeclareMathOperator{\Eff}{Eff}
 
 \DeclareMathOperator{\SSpec}{\cSpec}
 \DeclareMathOperator{\chara}{char}
 \DeclareMathOperator{\Divi}{Div}
 \DeclareMathOperator{\WDivi}{WDiv}
 \DeclareMathOperator{\ex}{ex}
 \DeclareMathOperator{\Ex}{Ex}
 \DeclareMathOperator{\Fan}{Fan}
 \DeclareMathOperator{\st}{st}
 \DeclareMathOperator{\cl}{cl}
 
 \DeclareMathOperator{\relint}{relint}
 \DeclareMathOperator{\Gal}{Gal}
 \DeclareMathOperator{\semistable}{ss}
 \DeclareMathOperator{\strictlysemistable}{sss}
 \DeclareMathOperator{\gp}{gp}
 \DeclareMathOperator{\free}{free}

\begin{document}
 
\maketitle
 
\begin{abstract}
The purpose of this paper is to study the geometry of images of morphisms from Mori dream spaces.  
First we prove that a variety which admits a surjective morphism from a Mori dream space
is again a Mori dream space.
Secondly we introduce a natural fan structure on the effective cone of Mori dream spaces.
We show that it encodes the information of Zariski decompositions, which
in turn is equivalent to the information of the variation of GIT quotients of their Cox rings.
Finally we show that under a surjective morphism between Mori dream spaces,
the fan of the target space coincides with the restriction of the fan of the source.
\end{abstract}

\tableofcontents

%

\section{Introduction}
The definition and a number of basic properties of
Mori dream spaces were established in \cite{MR1786494}.
Mori dream spaces are varieties whose line bundles
satisfy a strong condition. As proven in the same paper,
they are also characterized by the fact that their birational geometry
is described by the Variation of Geometric Invariant Theory
(VGIT for short) of its Cox ring.
Both of these two aspects of Mori dream spaces are of great use
and enrich their geometry.

Many interesting varieties are known to be Mori dream spaces.
Among them are toric varieties (and more generally spherical varieties \cite{MR1322696}),
varieties of Fano type in characteristic zero \cite[Corollary 1.3.2]{MR2601039}, and
K3 surfaces whose automorphism groups are finite \cite{MR2660680}.
More generally, the cone conjecture asserts that Calabi-Yau varieties have similar structures
as Mori dream spaces "up to the action of birational automorphism groups" \cite{MR1468356}.
In addition, several moduli spaces are known to be Mori dream spaces.
In nice situations the birational models turn out to be different modular compactifications,
and the space of parameters (more specifically that of stability conditions)
for the moduli problem is directly related to the Picard group.
A typical example of such is the moduli space of rank two parabolic bundles
on a pointed projective line \cite{mukai2005finite}.

One of the main problems concerning Mori dream spaces is to find new examples of them.
What makes it difficult is the subtlety of the property of being a Mori dream space.
In fact, it is not preserved by most of the natural operations:
see \pref{eg:failure_of_being_MDS_under_operations} for a list of such phenomena.

The purpose of this paper is to study the geometry of the images of morphisms from Mori dream spaces
and clarify the relationship between the geometry of the target space and that of the source.
Our first result is the following

\begin{theorem}\label{th:main}
Let $X,Y$ be normal $\bQ$-factorial projective varieties, and
$
f
\colon
X
\to
Y
$
be a surjective morphism.
If $X$ is a Mori dream space, then so is $Y$.
\end{theorem}

It was proven in \cite[Corollary 5.2]{MR2944479} that
the image of a variety of Fano type 
again is of Fano type.
Since varieties of Fano type (in characteristic zero) are Mori dream spaces,
\pref{th:main} can be regarded as a partial generalization of the result.
In fact another proof to \cite[Corollary 5.2]{MR2944479} was given in
\cite[Corollary 5.5]{MR3275656} by using \pref{th:main}.
We also expect \pref{th:main} will be useful to construct new examples of Mori dream spaces.

We prove \pref{th:main} by showing the finite generation of Cox rings of $Y$ over the base field,
using \pref{th:characterization} which says that the finite generation of Cox rings characterizes Mori dream spaces.
We deduce it from the finite generation of Cox rings of $X$ as follows.

By taking the Stein factorization of
$f$, the proof is divided into two parts: the case when $f$ is an algebraic fiber space
(\pref{sc:Multi-section_ring_and_Cox_ring})
and the case when $f$ is a finite morphism
(\pref{sc:Finite_generation_of_multi-section_rings_under_finite_morphisms}).
When $f$ is finite, we further decompose it into the separable part and the purely inseparable part,
and treat them independently (but with a somewhat similar idea). Combining them, we conclude the
proof of \pref{th:main} in \pref{sc:Proof of Theorem main}.

Next we compare the geometry of $Y$ with that of $X$.
For that purpose we introduce a fan structure on the effective cone of Mori dream space
(see \pref{dp:def of fan}).
We denote it by $\Fan{(\cdot)}$, and give two interpretations to it.

One is from the point of view of Zariski decompositions of line bundles.
As we will see in \pref{pr:ZD},
any pseudo-effective line bundle on a Mori dream space admits a Zariski decomposition
in the strongest sense.
We say that two line bundles on a Mori dream space are \emph{strongly Mori equivalent}
if the negative parts of their Zariski
decompositions have the same support and the positive parts define the same Iitaka fibration
(see \pref{df:strong Mori equivalence}). Then we prove that
two line bundles on a Mori dream space are strongly Mori equivalent if and only if they are contained in the
relative interior of the same cone of the fan (\pref{sc:Zariski decompositions and the fan}).

The other interpretation comes from the VGIT of Cox rings.
Cox rings are canonically graded by the divisor class group
(modulo some minor ambiguities),
so that the affine variety defined by the Cox ring comes with the natural
dual torus action of the divisor class group.
Since characters of the dual torus canonically correspond to divisor classes,
we can associate to a line bundle on a Mori dream space a character of the dual torus (see
\pref{sc:Strong Mori equivalence is equivalent to GIT equivalence}).
Thus we obtain a GIT problem.
We say that two line bundles on a Mori dream space are GIT equivalent if
the semi-stable loci of the corresponding characters coincide (see \pref{sc:Mori dream space and GIT revisited}).
Then we show in \pref{sc:Strong Mori equivalence is equivalent to GIT equivalence}
that two line bundles are GIT equivalent
if and only if they are strongly Mori equivalent.

Now let $f\colon X\to Y$ be a surjective morphism between Mori dream spaces.
By the interpretations of the fan above, the relation of the geometry of $Y$ with that of $X$ can be encoded as
the relation of $\Fan{(Y)}$ with $\Fan{(X)}$.
In order to see this, note that we can
regard $\Pic{(Y)}_{\bR}$ as a subspace of $\Pic{(X)}_{\bR}$ via the injection defined by the pull-back
\begin{equation*}
f^{*}\colon\Pic{(Y)}_{\bR}\hookrightarrow\Pic{(X)}_{\bR}.
\end{equation*}
Using this we can restrict $\Fan{(X)}$ to $\Eff{(Y)}$ and thus obtain the fan $\Fan{(X)}|_{\Eff{(Y)}}$.

\begin{theorem}\label{th:fans}
With the same assumptions as in \pref{th:main},
\begin{equation*}
\Fan{(Y)}=\Fan{(X)}\vert_{\Eff{(Y)}}
\end{equation*}
holds.
\end{theorem}
See \pref{eg:bl-up_of_P^3} for an illustration of \pref{th:fans}.
We give two proofs to \pref{th:fans} corresponding to the two characterizations of the
relative interiors of the cones explained above. These are treated respectively in
\pref{sc:Comparison of the fans without GIT} and
\pref{sc:Comparison of the fans via GIT}.
In both of the proofs \pref{th:finite} plays a key role.

In the final section, we extend our results to not necessarily $\bQ$-factorial Mori dream spaces and Mori dream regions.

Here we give several comments on relevant works.
In \cite{MR2811268}, it was shown that
the projective GIT quotient of an invariant open subset of a Mori dream space
by an action of a reductive group is a Mori dream space.
We quoted some ideas from the paper, and also there are overlapping results.

The idea of introducing fan structures on cones of divisors has a long history.
We can go back at least to \cite{Kawamata_CBU}, and later \cite{MR1420223}.
Our fan structure on the effective cone, restricted to the movable cone, coincides with the one
introduced in \cite{MR1786494}.
In \cite{MR1786494} they only consider the interiors of the maximal dimensional cones,
but the basic ideas needed for the proof of our refined results implicitly appear there.
For toric varieties, our fan structure was classically well known as the Gelfand-Kapranov-Zelevinsky
decomposition introduced by Oda and Park \cite{MR1117211}.

After the author finished the first version of the draft, he found that
J\"{u}rgen Hausen introduced in \cite{MR2499353} the notion of
GIT fan, starting from the VGIT of Cox rings. This seems to coincide with
our fan structure, which we define starting from the geometry of line bundles.

Finally, our treatment of Mori dream regions is not thorough.
After a draft of this paper appeared on arXiv, the paper
\cite{Kaloghiros:2012if} appeared and they gave a very systematic treatise.
The author believes that many part of their story can be
interpreted via VGIT, as we did in this paper.

\subsection*{Acknowledgements}
The author would like to express his gratitude to his advisor Yujiro Kawamata
for many suggestions, especially for asking him about the positive characteristic case.
The author is indebted to Tadakazu Sawada for informing him of basics and
examples of quotients by rational vector fields, and
to Young-Hoon Kiem for his insightful question.
He would like to thank Yoshinori Gongyo, Atsushi Ito, and Akiyoshi Sannai for
stimulating discussions. He would also like to thank the referees of this paper
for very carefully reading the manuscript and providing him with many suggestions and corrections.
The author was partially supported by Grant-in-Aid for JSPS fellows 22-849.

\subsection*{Conventions and Notations}
We work over a field $\bfk$. 
Unless otherwise stated, every variety is assumed to be projective,
geometrically integral and geometrically normal over $\bfk$.
For an abelian group $\Gamma$ and a field $K$, we denote by
$\Gamma_K$ the $K$-vector space $\Gamma\otimes_{\bZ}K$.
For notations and terminologies of Mori dream spaces and Cox rings we follow \cite{MR1786494}, and
for those of (V)GIT we follow \cite{Dolgachev-Hu} and \cite{Mumford-Fogarty-Kirwan}.


\section{Preliminaries on Mori dream space and the multi-section ring}
\label{sc:Preliminaries on Mori dream space and the multi-section ring}
\subsection{Mori dream space}\label{sc:Mori dream space}
We briefly recall definitions and some of the basic
notions/ properties about Mori dream spaces which we need in this paper.
See \cite{MR1786494} for details.

\begin{definition}\label{df:definition of cones}
Let $X$ be a projective variety over $\bfk$.
We denote by $\Num{(X)}$ the group of Cartier divisors on $X$ modulo
numerical equivalence. The cone spanned by nef divisors in $\Num{(X)}_{\bR}$ is denoted by
$\Nef{(X)}$.
Similarly, the closure of the cone spanned by movable
(resp. effective)
divisors is denoted by $\Eff{(X)}$ (resp. $\Mov{(X)}$).
\end{definition}

The closure operations in the definition of the cones $\Mov{(X)}$ and $\Eff{(X)}$
are unnecessary for Mori dream spaces, due to \pref{pr:[Proposition1.11(2)]{MR1786494}} below.
This is not the case for general varieties.

\begin{definition}\label{df:sqm}
Let $X$ be a normal projective variety. A \emph{small $\bQ$-factorial modification}
of $X$ is a small (i.e. isomorphic in codimension one) birational map
$f\colon X\dasharrow Y$ to another normal $\bQ$-factorial projective variety $Y$.
\end{definition}


\begin{definition}\label{df:Mori dream space} 
A normal projective variety $X$ is called  a \emph{Mori dream space}
provided that the following conditions hold.
\begin{enumerate}[(1)]
\item $X$ is $\bQ$-factorial and $\PPic^0_X$ has dimension zero.  \label{it:fin_pic}
\item $\Nef{(X)}$ is the affine hull of finitely many semi-ample
line bundles. \label{it:nef implies sa}
\item There is a finite collection of small $\bQ$-factorial modifications $f_i\colon X \dasharrow X_i$
such that each $X_i$ satisfies (\pref{it:fin_pic})(\pref{it:nef implies sa}) and $\Mov{(X)}$ is the union
of the $f_i^*\Nef{(X_i)}$.\label{it:movable cone condition}
\end{enumerate}
\end{definition}

\begin{remark}
The symbol $\PPic_X^0$ in (\pref{it:fin_pic}) denotes the identity component of
the Picard scheme of $X$ over $\bfk$
(see \cite[Chapter 9]{MR2222646} for details).
It has the following properties.
\begin{itemize}
\item It represents the relative Picard sheaf on the big \'{e}tale site $\text{\'{E}t}_{\bfk}$.

\item We always have an inequality
\begin{align*}
\dim H^1(X, \cO_X) \ge \dim \PPic^0_X.
\end{align*}
When
$\bfk$
is of characteristic zero,
the equality holds.
\end{itemize}

The condition (\pref{it:fin_pic}) above is slightly stronger than that of \cite[1.10. Definition]{MR1786494}. The former is preserved under change of base field,
whereas the latter is not.
Indeed, when
$
\bfk
=
\overline{\bfk}
$
and is not the closure of a finite field, then
these two conditions are equivalent by \cite[Theorem 10.1]{MR0337997}.
On the contrary, when
$
\bfk = \bar{\bF}_p
$
and $X$ is smooth over it, the map
\begin{align*}
\Pic{(X)}_{\bQ}\to \Num{(X)}_{\bQ}
\end{align*}
is \emph{always} an isomorphism.
Since it is natural to exclude curves of positive genera and abelian varieties
from Mori dream spaces, our condition (\pref{it:fin_pic}) seems to be
more appropriate than the original one.
\end{remark}

Let
$
D
$
be a
$
\bQ
$-Cartier divisor on a normal projective variety
$
X
$
whose section ring is finitely generated.
Then there exists a canonical dominant rational map
\[
\varphi_D
\colon
X
\dasharrow
\Proj R ( X , \cO_X ( D ) )
\]
such that the rational map
$
\Phi_{|D|}
\colon
X
\dasharrow
\bP H^0 ( X , \cO_X (D) )
$
factorizes through
$
\varphi_D
$.
If we replace
$
D
$
with
$
mD
$
for any
$
m
>0
$,
we see that there exists a canonical isomorphism
\[
\Proj R ( X , \cO_X ( D ) )
\simto
\Proj R ( X , \cO_X ( mD ) )
\]
which commutes with the rational maps
$
\varphi_D
$
and
$
\varphi_{mD}
$.
Combined with standard arguments for Iitaka fibrations
\cite[Sections 2.1B and 2.1C]{MR2095471},
this implies that
$
\Proj R ( X , \cO_X ( D ) )
$
has normal singularities and
$
\varphi_D
$
is an algebraic fiber space (in fact the Iitaka fibration)
in the sense that the subfield
\[
\varphi_D^*
\colon
\bfk ( \Proj R ( X , \cO_X ( D ) ) )
\hookrightarrow
\bfk ( X )
\]
is algebraically closed in
$
\bfk(X)
$.

\begin{definition}[{$=$\cite[Definitions 1.3 and 1.4]{MR1786494}}]
Let $D_1$ and $D_2$ be
two $\bQ$-Cartier divisors on $X$ with finitely generated
section rings.
We say $D_1$
and $D_2$ are \emph{Mori equivalent} if the dominant rational maps
$
\varphi_{D_1}
$
and
$
\varphi_{D_2}
$
are isomorphic, in the sense that there is an isomorphism
between their target spaces which makes the obvious triangular diagram commutative.
\end{definition}

\begin{definition}\label{df:Mori chamber}
Let $X$ be a normal projective variety satisfying \pref{df:Mori dream space} (\pref{it:fin_pic}).
A \emph{Mori chamber} of $X$ is the closure of a Mori equivalence class in $\Pic{(X)}_{\bR}$
with non-empty interior.
\end{definition}

In order to describe the Mori chambers of Mori dream spaces,
we recall some basic facts about birational contractions ($=$
birational maps surjective in codimension one).
Let
$
 g \colon X \dasharrow Y
$
be a birational contraction between normal projective varieties,
with
$X$
$\bQ$-factorial.
We denote by
\emph{$\ex{(g)}$}
the subcone of
$\Eff{(X)}$
spanned by the $g$-exceptional effective divisors.
Since any integral effective divisor $E\in\ex{(g)}$ satisfies $h^0(X,\cO_X(E))=1$,
we have the following easy fact:

\begin{lemma}\label{lm:negative cone}
The cone $\ex{(g)}$ is the simplicial cone
whose extremal rays are spanned by exceptional prime divisors of $g$.
In particular
$N_1 , N_2 \in \ex{(g)}$
have the same support if and only if
they are contained in the relative interior of the same face of $\ex{(g)}$.
\end{lemma}

The effective cone of a Mori dream space has a natural decomposition into Mori chambers
due to \cite[Proposition 1.11 (2)]{MR1786494}:

\begin{proposition}\label{pr:[Proposition1.11(2)]{MR1786494}}
Let
$X$
be a Mori dream space.
There are finitely many contracting birational maps $g_i\colon X\dasharrow Y_i$, with
$Y_i$ a Mori dream space, such that
\begin{equation*}
\Eff{(X)}=\bigcup_{i}g_i^{*}\Nef{(Y_i)}*\ex{(g_i)}
\end{equation*}
gives a decomposition of the effective cone into closed rational polyhedral subcones with
disjoint interiors. Each $g_i^{*}\Nef{(Y_i)}*\ex{(g_i)}$ is a Mori chamber of $X$.
\end{proposition}

\noindent
Here
$C*D$ denotes the join of the cones $C$ and $D$.

\subsection{The fan of a Mori dream space}\label{sc:The fan of a Mori dream space}

We introduce a fan structure on the effective cone of Mori dream space.
\begin{definition-proposition}\label{dp:def of fan}
Let $X$ be a Mori dream space. The set of faces of Mori chambers of $X$ forms a fan whose support coincides with
the effective cone of $X$. We denote it by $\Fan{(X)}$.
\end{definition-proposition}
\begin{remark}
The fan structure on $\Mov{(X)}$
introduced in \cite[Proposition 1.11(3)]{MR1786494}
is the restriction of $\Fan{(X)}$ to $\Mov{(X)}$.
\end{remark}
\begin{proof}
All we have to show is that the intersection of two cones of $\Fan{(X)}$ is a face of each cone.
By standard arguments of convex geometry,
it amounts to showing that for any Mori chambers $\cC_1$ and $\cC_2$, the intersection
$\cC_1 \cap \cC_2$
is a face of $\cC_2$.

Let $g_i\colon X\dasharrow Y_i$ ($i=1,2$) be the contracting birational map corresponding to $\cC_i$,
so that
$\cC_i=\cP_i*\cN_i$, where
$\cP_i=g_i^{*}\Nef{(Y_i)}$ and $\cN_i=\ex{(g_i)}$, holds.
The uniqueness of Zariski decomposition (see \pref{rm:uniqueness_of_ZD} below) implies the equality
\begin{align*}
\cC_1\cap\cC_2=(\cP_1\cap\cP_2)*
(\cN_1\cap\cN_2).
\end{align*}

As stated in \cite[Proposition 1.11(3)]{MR1786494},
we know that
$\cP_1\cap\cP_2$ is a face of $\cP_2$.
Similarly we can check that
$\cN_1\cap\cN_2$ is a face of $\cN_2$.
In fact, let $A=\sum a_iE_i$ ($a_i\ge 0$) and $B=\sum b_iE_i$ ($b_i\ge 0$) be
two elements of $\cN_2$ such that $A+B\in\cN_1$.
Since $h^{0}(X,\cO(A+B))=1$, we see
$\Supp{(A+B)}\subset \Ex{(g_1)}$. Hence we see $\Supp{(A)}, \Supp{(B)}\subset \Ex{(g_1)}$,
which implies
$A, B \in \scN_1$.
Therefore we see
$\cC_1\cap\cC_2$
is a face of $\scP_2 * \scN_2 = \cC_2$,
concluding the proof.
\end{proof}


\subsection{Zariski decompositions and the fan}\label{sc:Zariski decompositions and the fan}

Next we give an explicit description of Zariski decomposition
(which is slightly stronger than the sense of Cutkosky-Kawamata-Moriwaki)
of line bundles on a Mori dream space.
The existence of such a nice decomposition more or less
characterizes Mori dream spaces\footnote{the author would like to thank Professor
Y.-H. Kiem for asking him if it could be the case.}.


\begin{definition}\label{df:ZD}
Let
$X$
be a normal projective variety and
$D$
a pseudo-effective
$\bQ$-Cartier divisor on
$X$.
A \emph{Zariski decomposition} of the divisor $D$ is
a pair of
$\bQ$-Cartier divisors
$P$ and
$N$ on
$X$
satisfying the following conditions:
\begin{itemize}
\item
$P$ is nef.

\item
$N$ is effective.

\item
$D$
is
$\bQ$-linearly equivalent to
$P+N$.

\item
For any sufficient divisible
$m
\in
\bZ_{>0}$
the natural map
\begin{align}\label{eq:ZD}
\otimes s_{mN} \colon H^0(X, \cO_X(mP)) \to
H^0(X, \cO_X(mD)),
\end{align}
where $s_{mN}$ the tautological section of
$\cO_X(mN)$, is an isomorphism.
\end{itemize}
More generally if
$
X
$
is a Mori dream space,
we say that
$
D = P + N
$
is a Zariski decomposition of
$
D
$
if there exists a small
$
\bQ
$-factorial modification
$
\varphi
\colon
X
\dasharrow
X'
$
such that
$
\varphi_*D = \varphi_*P + \varphi_*N
$
is a Zariski decomposition in the above sense.
\end{definition}

\begin{remark}\label{rm:uniqueness_of_ZD}
On Mori dream spaces,
Zariski decomposition is unique.
In fact, by taking sufficiently divisible
$m>0$,
we can characterize the negative part
$N$
as
$\frac{1}{m}$
of the fixed part of the complete linear system
$|mD|$.
\end{remark}

\begin{proposition}\label{pr:ZD}
Let $X$ be a Mori dream space. Consider the decomposition of $\Eff{(X)}$ into
the Mori chambers given in \pref{pr:[Proposition1.11(2)]{MR1786494}}:
\begin{equation*}
\Eff{(X)}=\bigcup_{\textrm{finite}}\cC.
\end{equation*}
Then for each chamber $\cC$ there exists a small $\bQ$-factorial modification
$f_i\colon X\dasharrow X_i$ of $X$ and two $\bQ$-linear maps
\begin{equation*}
P,N\colon \cC\to \Eff{(X)}
\end{equation*}
such that for any $\bZ$-divisor $D\in \cC $,
$D\sim_{\bQ}P(D)+N(D)$ gives a Zariski decomposition of $D$ as a
divisor on $X_i$.

Conversely let $X$ be a normal projective variety satisfying \pref{df:Mori dream space} (\pref{it:fin_pic}).
Assume that
$\Eff{(X)}$
is decomposed into finitely many chambers $\cC$ on each of which
there exists $\bQ$-linear Zariski decompositions in the sense above, with
positive parts semi-ample on suitable small
$\bQ$-factorial modifications. Then $X$ is a Mori dream space.
\end{proposition}

\begin{proof}
Let $\cC$ be a Mori chamber, and
$g \colon X \dasharrow Y$
be the corresponding contracting birational map to another Mori dream space $Y$ as in
\pref{pr:[Proposition1.11(2)]{MR1786494}}.
We can replace $X$ with one of its small $\bQ$-factorial modifications so that $g$ becomes a morphism
by (\pref{it:movable cone condition}) of \pref{df:Mori dream space}.
Now we define the maps $P,N$ as follows:
\begin{itemize}
\item $P(D)=g^* g_*D$.
\item $N(D)=D-P(D)$.
\end{itemize}
Note that
$N(D)$
is a
$g$-exceptional effective
$\bQ$-divisor.
Since $h^{0}(X,\cO_X(mN(D)))=1$
holds for any sufficiently divisible
positive integer $m$,
the map
\pref{eq:ZD}
is uniquely defined up to constant.
When $m$ is sufficiently divisible so that $mP(D)$ is a $\bZ$-divisor,
it is easy to see that this map has the required properties.

The converse can be shown by checking the finite generation of Cox rings via exactly the
same arguments as in the proof of \pref{lm:afs}.
\end{proof}

We introduce a stronger version of the Mori equivalence relation, which is closely related to
the fan of Mori dream spaces defined above:
\begin{definition}\label{df:strong Mori equivalence}
Let $X$ be a Mori dream space.
Two line bundles $L$ and $M$ are said to be \emph{strongly Mori equivalent}
if they are Mori equivalent and
$
 \Supp{(N(L))} = \Supp{(N(M))}
$.
\end{definition}

Now we state the relationship between the notion of strong Mori equivalence and the
fan structure of Mori dream spaces.

\begin{proposition}\label{pr:fan vs strong Mori equivalence}
Let $L$ and $M$ be Mori equivalent line bundles on a Mori dream space $X$.
Then they are strongly Mori equivalent if and only if
their stable base loci coincide.
On any Mori dream space $X$, strong Mori equivalence classes coincide with the relative interiors of the cones of
$\Fan{(X)}$ and vice versa.
\end{proposition}

\noindent
See \cite[Definition 2.1.20]{MR2095471} for the definition of stable base loci.

\begin{proof}
`if' direction of the first statement is trivial, so we prove the `only if' direction.
Since $L$ and $M$ are Mori equivalent, the positive parts $P(L)$ and $P(M)$,
which are movable on $X$, are the pull-backs of some ample divisors under the same
contracting rational map $\varphi\colon X\dasharrow Z$.
Therefore the stable base loci of $P(L)$ and $P(M)$ are the same as the
locus of indeterminacy of the rational map $\varphi$.
Now note that the stable base locus of $L$ is the union of the support of $N(L)$ and
the stable base locus of $P(L)$.
Since the same thing holds for $M$, we obtain the conclusion.

Next we show the second statement.
Let $\cC$ be a Mori chamber and
$\cC=\cP*\cN$
be the Zariski decomposition of the chamber.
By an elementary fact on convex cones, the join of a face of $\cP$
with a face of $\cN$ is a face of $\cC$, and any face of $\cC$ is of this form.
Moreover if $C$ is a face of $\cC$ and $C=P * N$ is the decomposition,
we have the equality
\begin{equation*}
C^{\relint}=P^{\relint}*N^{\relint}.
\end{equation*}

Recall also that the relative interior of a face of $P$ is a strong Mori equivalence class.
This follows from the fact that two semi-ample line bundles are Mori equivalent if and only if
the set of curves contracted by the morphisms coincide.
Also the same thing hold for faces of $N$ by \pref{lm:negative cone}.
In particular, line bundles contained in 
$C^{\relint}$
are strongly Mori equivalent to each other.

Conversely suppose that two line bundles $L, M$ are strongly Mori equivalent.
Let $g \colon X\dasharrow Y$ be a birational contraction whose Mori chamber
contains $L$, so that the rational map
$\varphi_L=\varphi_M \colon X \dasharrow Z$
factors through a morphism $\psi\colon Y\to Z$.
Since the positive part $P(M)$ is the pull-back of an ample divisor on $Z$,
it is the pull-back of a semi-ample divisor on $Y$. Together with the equality
$\Supp N(L)=\Supp N(M)$, this implies that $M$ is also contained in the Mori chamber of $g$.
The rest of the proof follows from the arguments above.
\end{proof}


\subsection{Multi-section ring and the Cox ring}
\label{sc:Multi-section_ring_and_Cox_ring}

The notion of multi-section ring is of fundamental importance in the study of Mori dream spaces.

\begin{definition}\label{df:multi-section ring}
Let $X$ be a normal variety with
$H^0(X,\cO_{X})=\bfk$.
Let $\Gamma\subset\WDivi{(X)}$ be a semigroup of Weil divisors.
The \emph{multi-section ring} $R(X, \Gamma)$ \emph{associated to} $\Gamma$ is
the $\Gamma$-graded $\bfk$-algebra defined by
\begin{equation*}
R(X, \Gamma)=\bigoplus_{D\in\Gamma}H^0(X,\cO_X(D)).
\end{equation*}
In particular for a divisor $D$ on $X$, we define the section ring of $D$ by
\begin{equation*}
 R \lb X, \cO _X ( D ) \rb 
 = R \lb X, \bZ _{ \ge 0 } D  \rb
 = \bigoplus _{ m \ge 0 } H ^0  \lb X, \cO _X ( m D ) \rb.
\end{equation*}
\end{definition}

We define the notion of Cox rings as multi-section rings.

\begin{definition}\label{df:Cox ring}
Let $X$ be a normal projective variety such that the Weil divisor class group $\Cl{(X)}$
is finitely generated.
A \emph{Cox ring} of $X$ is the multi-section ring
\begin{equation*}
R(X, \Gamma)=\bigoplus_{D\in\Gamma}H^0(X,\cO_{X}(D)),
\end{equation*}
where $\Gamma\subset\WDivi{(X)}$ is a finitely generated group of Weil divisors such that
the class map
$\Gamma_{\bQ} \to \Cl{(X)}_{\bQ}$
is an isomorphism.
\end{definition}

\begin{remark}
Although Cox rings defined as above depend on the choice of the group $\Gamma$,
basic properties such as finite generation is independent of the choice.

When the divisor class group is torsion free, it is standard to choose such a
$\Gamma$
that maps onto the divisor class group.
For those
$\Gamma$,
the resulting Cox rings are all isomorphic to each other.
Even if the divisor class group has torsion,
in fact there still exists a sophisticated way to define \emph{the} Cox ring,
which is graded by the divisor class group and
is unique up to isomorphism (see \cite[Section 1.4.2]{MR3307753}).
Our Cox rings
$ R( X, \Gamma)$
are subrings of that Cox ring of \cite{MR3307753},
and the finite generation of Cox rings in our sense is equivalent to
that of theirs.
We adopt our definition in this paper since
it is simpler to handle and enough for our purpose.
See \cite[Remark 2.19]{MR3275656} for more details.
\end{remark}

The birational geometry of Mori dream spaces are tautologically equivalent
to the VGIT of Cox rings. This was first proven in \cite{MR1786494}, and
the equivalence will be treated in detail in \pref{sc:Mori dream space and GIT revisited}.
The most important consequence is the following characterization of
Mori dream spaces via the finite generation of their Cox rings.

\begin{theorem}[{$=$ \cite[Proposition 2.9]{MR1786494}}]\label{th:characterization}
A normal projective variety satisfying
\pref{df:Mori dream space} (\pref{it:fin_pic}) is a Mori dream space if and only if its Cox ring
is of finite type over $\bfk$.
\end{theorem}

Here we prove the finite generation of multi-section rings on Mori dream spaces,
which is a generalization of the `if' direction of \pref{th:characterization}.
This was first proven by B\"aker \cite[Theorem 1.2]{MR2811268} by using the finite generation theorem for
invariant subrings for the actions of reductive groups. We give an alternative geometric
proof based on Zariski decompositions on Mori dream spaces.

\begin{lemma}\label{lm:afs}
Let $X$ be a Mori dream space. Let $\Gamma\subset\WDivi{(X)}$ be a finitely generated group of Weil divisors.
Then the multi-section ring $R(X, \Gamma)$ is of finite type over $\bfk$.
\end{lemma}
\begin{proof}
Without loss of generality, we will freely replace $\Gamma$ with a subgroup $\Gamma'$ of finite index.
In fact, since $R(X, \Gamma')\subset R(X, \Gamma)$ is an integral extension,
the finite generation of $R(X, \Gamma')$ implies that of $R(X, \Gamma)$
(see \cite[Corollary 13.13]{Eisenbud_CA}).

Consider the natural map $\cl \colon \Gamma\to\Cl{(X)}$.
By replacing
$\Gamma$
with a subgroup of finite index if necessary,
we can assume that
$\Image{\cl}$
is a free abelian group of finite rank.
Thus we obtain a decomposition
$
 \Gamma
=
\ker\cl
\oplus
\Gamma_1
$
so that
$\cl|_{\Gamma_1}$
is an isomorphism onto its image.
Since we have an isomorphism
$
 R(X, \Gamma)
\simeq
R(X, \Gamma_1)[\ker\cl],
$
it is enough to prove the finite generation of
$R(X, \Gamma_1)$. Replace
$\Gamma$ with $\Gamma_1$,
so that we can assume
$\cl$ is injective and
$\Gamma$ is torsion free.

Let $\cC$ be a Mori chamber.
Since $\cC$ is a rational polyhedral cone, the semigroup
$
 \Gamma_{\cC}=\Gamma\cap \cC
$
is finitely generated.
Let
$g \colon X\dasharrow Y$ be the birational contraction corresponding to the Mori chamber $\cC$.
Recall from the proof of \pref{pr:ZD} that for any $D\in\Gamma_{\cC}$ we have
a Zariski decomposition
\begin{equation}\label{eq:geometric_ZD}
D=g^* g_* D+(D-g^* g_* D).
\end{equation}
Since $\Gamma$ is finitely generated and there are only finitely many Mori chambers,
there exists a positive integer $m>0$ such that
for any Mori chamber $\cC$ and any divisor $D\in (m\Gamma)_{\cC}$,
the positive and the negative parts of the decomposition \pref{eq:geometric_ZD}
are both $\bZ$-divisors.
As before, replace $\Gamma$ with $m\Gamma$.

With these preparations, we have an isomorphism
\begin{equation*}
 \varphi \colon R ( X, \Gamma _{ \cC } ) \simto
 R (Y, P(\Gamma_{\cC}))[N(\Gamma_{\cC})]
\end{equation*}
which divides a section of the divisor
$D$
into the positive and the negative part according to
\pref{eq:geometric_ZD}.
Because of the existence and the uniqueness of Zariski decomposition,
$\varphi$ is an isomorphism.

Now since $P(\Gamma_\cC)$ is a finitely generated semigroup of semi-ample divisors
on $Y$ and 
$N(\Gamma_\cC)$ is a finitely generated semigroup, we see that $R(X, \Gamma_\cC)$ is
of finite type over $\bfk$ (see \cite[Lemma 2.8]{MR1786494}).
Since there are only finitely many chambers
$\cC$,
$R(X, \Gamma)$ itself is finitely generated.
\end{proof}


\section{Finite generation of multi-section rings under finite morphisms}
\label{sc:Finite_generation_of_multi-section_rings_under_finite_morphisms}
In this section we prove that finite generation of multi-section rings is preserved
under finite morphisms.

\begin{theorem}\label{th:finite}
Let $f\colon X\to Y$ be a finite surjective morphism between normal varieties, and
$\Gamma\subset \WDivi{(Y)}$ a finitely generated semigroup of Weil divisors.
Then $R(Y, \Gamma)$ is of finite type over $\bfk$ if and only if $R(X, f^{*}\Gamma)$ is.
Moreover, in this situation,
the natural extension of multi-section rings $R(Y, \Gamma)\subset R(X, f^{*}\Gamma)$ is finite.
\end{theorem}


\subsection{Preliminary for the proof of \pref{th:finite}}
\label{sc:Preliminary for the proof of Theorem finite}

In the proof of \pref{th:finite}, we frequently use the notion of universal torsors.
\begin{definition}\label{df:universal torsor}
Let $\Gamma\subset \WDivi{(Y)}$ be a finitely generated semigroup of Weil divisors.
We define \emph{the universal torsor associated to $\Gamma$} as
\begin{equation*}
\cS_Y(\Gamma)=\bigoplus_{D\in\Gamma}\cO_Y(D).
\end{equation*}
\end{definition}

\begin{remark}\label{rm:remarks_on_universal_torsor}
We use the following basic properties frequently.
 
\begin{enumerate}[(1)]
\item\label{it:R_X and S_X}
$H^0(Y,\cS_Y(\Gamma))=R(Y, \Gamma)$.

\item\label{it:flat} If $Y$ is non-singular, $\cS_Y$ is a smooth $\cO_Y$-algebra.

\item\label{it:S_X}
Given a morphism
$
 f \colon X \to Y
$,
we have
$
 R(X, f^{*} \Gamma ) = H^0 ( Y, f_{*} \cS_X ( f^{*} \Gamma ) )
$.
By the projection formula, we also obtain the isomorphism
$
 f_{*}\cS_X(f^{*}\Gamma)\simeq\cS_Y(\Gamma)\otimes_{\cO_Y}f_{*}\cO_X
$.

\item\label{it:removable}
Since $Y$ is normal, $R(Y, \Gamma)$ is unchanged if we replace $Y$ with its non-singular locus.
Hence when $f\colon X\to Y$ is a finite surjective morphism, by removing suitable
closed subsets of codimensions at least two, we can assume that $X,Y$ are non-singular.
\end{enumerate}
\end{remark}

\vspace{5mm}

Now we go back to the proof \pref{th:finite}.
Consider the factorization
$
 X \xxto[]{g} \Ytilde \xxto[]{h} Y
$
of $f$ into the purely inseparable part $g$
and the separable part $h$
(recall that $\Ytilde$ is obtained as the normalization of $Y$
in the separable closure of $k(Y)$ in $k(X)$).
We will treat $g$ and $h$ separately in the following two subsections,
although the idea of the arguments are basically the same.


\subsection{Purely inseparable case}\label{sc:Purely inseparable case}

Assume that $f$ is a purely inseparable finite morphism. 
Throughout this section we take the base change from
$\bfk$
to the algebraic closure
$\overline{\bfk}$ so that we can assume
the base field is algebraically closed.
Since we are only interested in
the finite generation of algebras of the form
$R( X , \Gamma)$ over $\bfk$,
this does not affect the arguments.

We divide the extension $\bfk(Y)\subset \bfk(X)$ into
subextensions of degree $p$, so that we can assume $\deg{(f)}=p$.
In order to study such morphisms, instead of the Galois correspondence
for Galois extensions,
we can use the following \emph{Jacobson correspondence}
(see \cite[Section 10]{MR3114931} for details).
Recall that a purely inseparable finite morphism
$f \colon X \to Y$
of normal varieties is
\emph{of height one}
if the corresponding extension of function fields satisfy
$
 \bfk{(X)}^p \subset
\bfk{(Y)} \subset
\bfk{(X)}.
$

\begin{proposition}[{$=$ \cite[Theorem 10.5]{MR3114931}}]\label{pr:Jacobson correspondence}
Let
$X$
be a smooth variety over $\bfk$.
Then there is a canonical bijection between the following two sets.
\begin{enumerate}
\item Set of finite morphisms
$\varphi \colon X \to Y$
of height one.

\item Set of $p$-closed foliations in $T_X$.
\end{enumerate}
\end{proposition}

A foliation in $T_X$ is a saturated subsheaf which is, regarded as a sheaf of
derivations, closed under the bracket
$[ \  , \ ]$.
It is said to be $p$-closed if it is also closed under the
$p$-th power operation: i.e.,
$\delta \mapsto \delta^{[p]}:= \delta \circ \delta \circ \cdots \circ \delta$ ($p$ times).

\begin{corollary}\label{cr:insep_quot}
Let $f\colon X\to Y$ be a purely-inseparable finite morphism of degree
$p$ between normal varieties over $\bfk$, with $X$ smooth.
Then there exists a $p$-closed foliation
$\scE \subset T_X$ such that
\begin{equation*}
\cO_Y = \cO_X^{\scE}= \{x \in \cO_X | \scE{(x)}=0\}.
\end{equation*}

Moreover for any smooth morphism $Z\to Y$,
the base change morphism
$
 X \times_{Y} Z \to Z
$
is again a purely-inseparable morphism
of degree $p$
corresponding to the pull-back of
$\scE$.
\end{corollary}

\begin{proof}
Since
$f$
is purely inseparable of degree $p$, we easily see that it is of height one.
Thus the first statement follows from \pref{pr:Jacobson correspondence}.
The second statement follows from the functorial nature of the
Jacobson correspondence. See \cite[Proposition 2.3]{MR927978}.
\end{proof}

\begin{proof}[Proof of \pref{th:finite} when $f$ is purely inseparable]
Suppose that $f$ is purely inseparable. As mentioned before, we can assume $\deg{(f)}=p$
and
$\bfk$
is algebraically closed.
By \pref{rm:remarks_on_universal_torsor} (\pref{it:removable})
and (\pref{it:flat}),
we can assume that $X,Y$ are non-singular and hence
$
 Z:=\SSpec_Y\cS_Y(\Gamma)\to Y
$
is smooth.
Therefore we can apply
\pref{cr:insep_quot} to the morphisms
$ f \colon X \to Y$ and
$Z\to Y$, so that
$W = X \times_Y Z \to Z$
is the quotient by a $p$-closed foliation on $W$.
Therefore we obtain the inclusions
\begin{align}\label{eq:height_one_for_structure_sheaves}
\cO_W^p \subset \cO_Z \subset \cO_W.
\end{align}

Note that
$
 W \simeq \SSpec_X \scS_X{(f^*\Gamma)}
$
by \pref{rm:remarks_on_universal_torsor}
(\pref{it:S_X}),
so that we see
$H^0(Z,\cO_Z)=R(Y, \Gamma)$
and
$H^0(W,\cO_W)=R(X, f^{*}\Gamma)$
by \pref{rm:remarks_on_universal_torsor} (\pref{it:R_X and S_X}) and (\pref{it:S_X}).
Hence the inclusions \pref{eq:height_one_for_structure_sheaves}
imply the inclusions
$
 R( X , f^* \Gamma)^p \subset
R (Y , \Gamma) \subset
R( X , f^* \Gamma).
$
From this we easily obtain the conclusions.
\end{proof}


\subsection{Separable case}\label{sc:Separable case}

Assume that $f$ is separable.
This case is relatively easier; by
passing to a Galois closure, we can describe $Y$ as a uniform geometric
quotient by the Galois group, so that we can apply
the finite generation theorem for invariant subrings.

\begin{proof}[Proof of \pref{th:finite} when $f$ is separable]
Suppose that $f$ is separable.
Let $L$ be the Galois closure of the extension $\bfk(Y)\subset \bfk(X)$, and let $W$ be the
normalization of $X$ in $L$ so that $L=\bfk(W)$ holds.
Note that the finite map $W\to Y$ is the quotient of $W$ by the Galois group $\Gal{(W/Y)}$.
By removing suitable closed subsets, we can assume that $X,Y,W$ are all non-singular
(see \pref{rm:remarks_on_universal_torsor} (\pref{it:removable})).

Since $W/X$ also is Galois, it is the uniform geometric quotient of the action of $\Gal{(W/X)}$
on $W$ (see the proof of \cite[Lemma A.1.2]{Mumford-Fogarty-Kirwan} for the uniformity). Non-singularity of $X$ implies that
$\cS_X{(\Gamma)}$ is a flat $\cO_X$-algebra
and hence we have
$
 \cS_X{(\Gamma)}=\cS_W{(\Gamma)}^{\Gal(W/X)}.
$
In particular we see that
$R(X, \Gamma)=R(W, \Gamma)^{\Gal(W/X)}$,
which means that $R(W, \Gamma)$ is an integral extension of
$R(X, \Gamma)$. By the same argument, we can also show
$R(Y, \Gamma)=R(W, \Gamma)^{\Gal{(W/Y)}}$.

Now suppose that $R(X, \Gamma)$ is of finite type over $\bfk$.
By the finiteness theorem for integral closures \cite[Corollary 13.13]{Eisenbud_CA}, we see
$R(W, \Gamma)$ also is of finite type over $\bfk$.
This in turn implies the finite generation of $R(Y, \Gamma)$, since $\Gal{(W/Y)}$ is a finite group.
By similar arguments we can also check that the finite generation of $R(Y, \Gamma)$ implies that
of $R(X, \Gamma)$, concluding the proof of `if' part. `only if' part can be proven similarly, and
the finiteness of the extension $R(Y, \Gamma)\subset R(X, f^{*}\Gamma)$ follows also from
\cite[Corollary 13.13]{Eisenbud_CA}.
\end{proof}


\section{Proof of \pref{th:main}}\label{sc:Proof of Theorem main}

In this section we prove \pref{th:main} using \pref{th:characterization}.
First of all we check that $Y$ satisfies the condition (\pref{it:fin_pic}) of \pref{df:Mori dream space}.

\begin{lemma}\label{lm:pic}
Under the assumptions of \pref{th:main},
$\PPic_Y^0$ is also zero-dimensional.
\end{lemma}
\begin{proof}
Since we are only interested in dimensions, we
take the base change by $\bfk\subset\overline{\bfk}$ so that
we can assume the base field $\bfk$ is algebraically closed.
Let
$
 X\xxto[]{g}\Ytilde\xxto[]{h}Y
$
be the Stein factorization of $f$.
Since $f^*$ maps $\PPic^0_Y$ to $\PPic_X^0$ and the latter consists of a point
by the assumption, for any line bundle $L\in\Pic^0_Y(\bfk)$ we see $f^*L\simeq \cO_X$.
Since $g$ is an algebraic fiber space, this implies that $h^*L\simeq \cO_{\Ytilde}$.
Hence we see $L^{\otimes \deg h}\simeq\cO_Y$, and this implies that the morphism
$f^*\colon\PPic^0_Y\to\PPic_X^0$ is finite since both $\PPic^0_Y$ and $\PPic_X^0$
are projective schemes over $\bfk$. Therefore $\dim\PPic^0_Y\le\dim\PPic^0_X=0$,
concluding the proof.
\end{proof}

\begin{remark}\label{rm:direct}
For a surjective morphism $f\colon X\to Y$ between normal projective varieties,
$f^{*}\colon\Pic{(Y)}_{\bR}\to\Pic{(X)}_{\bR}$ is injective. We regard $\Pic{(Y)}_{\bR}$ as
a subspace of $\Pic{(X)}_{\bR}$ via the map $f^{*}$. Then from similar arguments 
as in the proof of \pref{lm:pic}, we can easily check the following statements
under the assumptions of \pref{th:main};
\begin{itemize}
\item
$\Eff{(Y)}=\Eff{(X)}\cap\Pic{(Y)}_{\bR}$,

\item
$\Nef{(Y)} = \Nef{(X)}\cap \Pic{(Y)}_{\bR}=\Nef{(X)}\cap\Eff{(Y)}.$
\end{itemize}
\end{remark}

\begin{proof}[Proof of \pref{th:main}]
By \pref{lm:pic} and \pref{th:characterization},
it is enough to show the finite generation of a Cox ring of $Y$.

Take a group $\Gamma\subset\Divi{(Y)}$ which maps isomorphically onto
a subgroup in $\Pic{(Y)}$ of finite index.
By \pref{lm:afs}, we know that $R(X, f^{*}\Gamma)$
is of finite type over $\bfk$. On the other hand
since $g$ is an algebraic fiber space, we have an isomorphism
$
 R(X, f^{*}\Gamma)\simeq R(\Ytilde, h^{*}\Gamma).
$
Since $h$ is finite, by applying \pref{th:finite},
we can conclude that $R(Y, \Gamma)$ is also of finite type over $\bfk$.
\end{proof}


\section{Comparison of the fans without GIT}\label{sc:Comparison of the fans without GIT}

In this section we prove \pref{th:fans} via direct arguments.
The problem is reduced to the following
\begin{theorem}\label{th:comparing strong Mori equivalence}
Let $f\colon X\to Y$ be a surjective morphism between Mori dream spaces.
Then two line bundles $L$ and $M$ on $Y$ are strongly Mori equivalent
if and only if $f^{*}L$ and $f^{*}M$ are
strongly Mori equivalent.
\end{theorem}

\noindent
See \pref{df:strong Mori equivalence} for the notion of strong Mori equivalence.
We first check that \pref{th:fans} actually follows from this.
\begin{proof}[Proof of \pref{th:fans}]
Take any $\sigma\in\Fan{(Y)}$. By \pref{pr:fan vs strong Mori equivalence} and \pref{th:comparing strong Mori equivalence},
there exists a cone $\Sigma\in\Fan{(X)}$ such that
$
 \sigma^{\relint}=\Sigma^{\relint}\cap\Eff{(Y)}
$.
Since the RHS is not empty, we can check
$
 \sigma
 =
 \overline{\lb \Sigma^{\relint}\cap\Eff{(Y)} \rb}=\Sigma\cap\Eff{(Y)}
$.

Conversely, let $\Sigma\in\Fan{(X)}$ be a cone which intersects $\Eff{(Y)}$.
Let $\Sigma'$ be the largest face of $\Sigma$ such that
$
 \Sigma\cap\Eff{(Y)}=\Sigma'\cap\Eff{(Y)}
$.
Note that
$
 \Sigma'^{\relint}\cap\Eff{(Y)}\not=\emptyset
$.
Again by \pref{pr:fan vs strong Mori equivalence} and \pref{th:comparing strong Mori equivalence},
there exists a cone $\sigma\in\Fan{(Y)}$ such that
$
 \Sigma'^{\relint}\cap\Eff{(Y)}=\sigma^{\relint}
$.
Taking the closures, we obtain
$
 \Sigma\cap\Eff{(Y)}=\Sigma'\cap\Eff{(Y)}=\sigma
$
and thus conclude the proof.
\end{proof}

To prove \pref{th:comparing strong Mori equivalence},
we begin with the following general result on stable base loci.

\begin{proposition}\label{pr:pull-back_and_stable_base_locus}
Let
$f \colon X \to Y$
be a surjective morphism between normal projective varieties.
For any Cartier divisor
$D$ on $Y$, we have the equality
$
 \bB{(f^* D)}
=
f^{-1} ( \bB(D) ).
$
\end{proposition}

\begin{proof}
By taking the Stein factorization,
we can assume
$f$ is either an algebraic fiber space or finite.
When
$f$ is an algebraic fiber space, the assertion is trivial
because of the isomorphism
$
 H^0 ( X , \cO_X( m f^* D ))
\simeq
H^0 ( Y, \cO_Y( m D )).
$

When
$f$ is finite, as we did in the proof of
\pref{th:finite},
we can further assume either
$f$ is purely inseparable of degree
$p$
or separable.
In the latter case, by \pref{cr:insep_quot},
we have the inclusions
$
 R ( X , D ) ^p
\subset
R ( Y , D )
\subset
R ( X , D ).
$
By standard arguments, the conclusion follows from this.

When
$f$
is separable,
by taking the Galois closure as we did in
\pref{sc:Separable case},
we can reduce it to the case when
$f$ is Galois.
This follows from the fact that
for a Galois morphism
$f \colon X \to Y$,
both
$\bB{(f^*D)}$
and
$f^{-1}(\bB{(D)})$
are $\Gal{(X/Y)}$-invariant.

Now let us assume $f$ is Galois with Galois group
$G$.
Assume that
$x \in X$ is not contained in
$\bB {(f^*D)}$, so that
there exists
$s \in H^0 ( X , \cO_X ( m f^* D))$
for some
positive integer
$m$
such that
$ s ( x ) \not= 0$.
Then there exists at least one homogeneous symmetric polynomial $p$
of the sections 
$g^* s$ \ ($g \in G$)
such that
$p(x) \not= 0$.
Since
$p \in H^0 (X, \cO_X(\deg p \cdot m f^*D))$
is
$G$-invariant,
it descends to a global section of
$\cO_Y(\deg p \cdot mD)$. Hence $x\not\in f^{-1}(\bB{(D)})$.
Since the inclusion $\bB {(f^*D)} \subset f^{-1}(\bB{(D)})$
is obvious, we conclude the proof.
\end{proof}

\begin{proof}[Proof of \pref{th:comparing strong Mori equivalence}]
By \pref{pr:pull-back_and_stable_base_locus},
in view of
\pref{pr:fan vs strong Mori equivalence},
it is enough to show the following claim.
\end{proof}

\begin{claim}\label{cl:comparing Mori equivalence}
Two line bundles
$L$ and
$M$ in
$\Eff{(Y)}$
are Mori equivalent if and only if $f^*L$ and $f^*M$ are
Mori equivalent.
\end{claim}

\begin{proof}
Taking the Stein factorization, we can assume that
$f$ is either an algebraic fiber space or finite.
Suppose $f$ is an algebraic fiber space and
$f^* L$
and
$f^* M$
are Mori equivalent line bundles.
Then we obtain the following commutative diagram.

\[
\xymatrix{
\Proj{R(X, f^{*}M)} \ar[rr]^{\sim} \ar[dd] & & \Proj{R(X, f^{*}L)} \ar[dd]\\
& X \ar@{-->}[lu]^{\varphi_{f^{*}M}} \ar@{-->}[ru]_{\varphi_{f^{*}L}} \ar[dd]^{f}& \\
\Proj{R(Y, M)}   & & \Proj{R(Y, L)} \\
& Y \ar@{-->}[lu]^{\varphi_{M}} \ar@{-->}[ru]_{\varphi_{L}} & \\}\]

\noindent
The top horizontal arrow in the diagram is an isomorphism which makes the
upper triangle commutative. 
Note that the two side vertical morphisms are isomorphisms, since $f$ is an
algebraic fiber space. Therefore $M$ and $L$ are Mori equivalent.
By similar arguments, we can also check that the equivalence of $L$ and $M$ implies that of
$f^*L$ and $f^*M$.

Now assume that $f$ is finite and $f^*L$ is equivalent to $f^*M$.
Consider the same diagram as above.
Since $R(X, f^*L)$ is finite over $R(Y, L)$ (resp. for $M$) by \pref{th:finite}
(consider
$\Gamma
=\bZ_{\ge 0}L$
so that
$R(Y, L)=R(Y, \Gamma)$),
we see that the two side vertical morphisms are also finite.

In order to prove the existence of an isomorphism from $\Proj{R(Y, M)}$ to $\Proj{R(Y, L)}$
which is compatible with any other maps,
it is enough to show that the morphism from
$
\Proj R ( X, f^*M )
$
to
$
\Proj R ( X , f^* L )
$
descends to a morphism from
$
\Proj R ( Y , M )
$
to
$
\Proj R ( Y , L )
$
and that the same thing holds if we interchange $M$ and $L$.
Recall the following decomposition of the morphism $f$ in
\pref{sc:Finite_generation_of_multi-section_rings_under_finite_morphisms}:

\[
\xymatrix{
 & W \ar[d] \ar[rd] &\\
X \ar[r] & S \ar[r] & Y\\}\].

\noindent
In the diagram above, $S$ is the separable closure of the extension $X/Y$ and $W$ is the Galois closure
of $S/Y$.

$\Proj{R(Y, L)}$ is obtained from $\Proj R(X, f^*L)$ as follows.
We repeatedly take quotients by
$p$-closed foliations which corresponds to a chain of degree $p$ subextensions of $X\to S$
(see \pref{cr:insep_quot}),
take its normalization in $\bfk(W)$ (see \cite[Example 2.1.12]{MR2095471}),
and take the quotient by the Galois group $G(W/Y)$
(see the arguments in \pref{sc:Separable case}).
Note that this process depends only on the initial data $\Proj R(X, f^*L)$,
subextensions of $X/Y$ and $W$. Therefore the isomorphism
between $R(X, f^{*}M)$ and $R(X, f^{*}L)$ descends to an isomorphism between
$\Proj{R(Y, M)}$ and $\Proj{R(Y, L)}$.

Conversely since $\Proj R(X, f^*L)$ is the normalization of $\Proj R(Y, L)$ in $\bfk(X)$,
the Mori equivalence of $L$ to $M$ implies that of $f^*L$ to $f^*M$. Thus we conclude the proof.
\end{proof}


\section{Mori dream space and GIT revisited}\label{sc:Mori dream space and GIT revisited}

In this section we explain the relations between the variation of GIT quotients (VGIT for short) of Cox rings
and the geometry of Mori dream spaces, to complement \cite{MR1786494}.

\begin{remark}
In \cite{MR1786494}, results on VGIT were quoted from \cite{Dolgachev-Hu}, 
in which varieties under group actions are assumed to be proper.
Since we have to deal with affine varieties defined
by Cox rings, we need VGIT for affine varieties.
One big difference from the proper case is the fact that
for a 1-parameter subgroup $\lambda$ of the group and a point $x$ on the variety,
the limit point $\lim_{t\to 0}\lambda (t)x$ may not exist.
As a consequence the wall defined by a point $x$ may not be a convex set, contrary to the proper case:
in fact, as we will see in the next subsection, the wall defined by
$x$ is the boundary of a certain convex polyhedral cone.
\end{remark}

\subsection{VGIT of torus actions on affine varieties}
\label{sc:VGIT of torus actions on affine varieties}
Let $G$ be a reductive group acting on a normal affine variety $V$.
Assume for simplicity that
only finitely many elements of $G$ acts on $V$ trivially.

Let $\chi{(G)}=\Hom_{\gp}(G,\bG_m)$ be the character group of $G$
and $\chi_{\bullet}{(G)}$ the group of 1-parameter subgroups of $G$,
so that we have the natural pairing
$
 \la \chi, \lambda\ra=n,
$
where $\chi\in\chi{(G)}$, $\lambda\in\chi_{\bullet}{(G)}$ and $(\chi\circ\lambda)(t)=t^{n}$.

For
$\chi\in\chi{(G)}$,
we denote by
$L_{\chi}$
the
$G$-linearized line bundle
$\cO_V$ with the twist by
$\chi$
along the fiber direction.
Let $U_{\chi}:=V^{\semistable}(L_{\chi})$
be the semi-stable locus of
$L_{\chi}$ on $V$,
and
$q_{\chi}\colon U_{\chi}\to Q_{\chi}=U_{\chi}//G$
be the quotient map.

\begin{proposition}[{$=$ \cite[Proposition 2.5]{King}}]\label{pr:King}
Let $G$ and $V$ as above, and
$\chi\in\chi (G)$ a character of $G$.
Then 
$x\in V$ is $L_{\chi}$(-semi)-stable if and only if
$\la \chi, \lambda\ra > 0 \ (\ge 0)$
holds for any 1-parameter subgroup
$\lambda\in\chi_{\bullet}(G)\setminus\{0\}$
such that $\lim_{t\to 0}\lambda (t)\cdot x$ exists.

\end{proposition}

Assume now that $G$ is an algebraic torus $T$.
Fix any $T$-equivariant embedding $V\subset \bA$ into
an affine space.
Let 
$\bA=\bigoplus_{\chi\in\chi{(T)}}\bA_{\chi}$
be the eigenspace decomposition.
Take a point $x\in V\subset\bA$, and the corresponding
decomposition $x=\sum_{\chi} x_{\chi}$ in $\bA$.
Since $V$ is closed in $\bA$, we see

$\begin{array}{ll}
& \lim_{t\to 0}\lambda (t)\cdot x \ \ \textrm{exists in V} \\
\iff & \lim_{t\to 0}\lambda (t)\cdot x \ \ \textrm{exists in $\bA$} \\
\iff & \la \chi, \lambda\ra\ge 0 \ \ \textrm{holds for all} \ \chi \ \textrm{such that} \ x_{\chi}\not=0. \\
\end{array}$

\noindent Consider the state set $\st{(x)}=\{\chi\in\chi{(T)}|x_{\chi}\not=0\}\subset\chi{(T)}$
of the point $x\in V$.
Denote by $\mathcal{D}_x\subset\chi{(T)}_{\bR}$ the cone spanned by $\st{(x)}$. Then
\begin{proposition} For a character $\chi\in\chi{(T)}$, 
$x\in V$ is $L_{\chi}$(-semi)-stable if and only if
$\chi\in\mathcal{D}_x^{\circ}$ (resp. $\scD_x$).
\end{proposition}

\begin{proof}
This is almost tautological. From \pref{pr:King} and the argument above,
$x\in V$ is 
$\chi$-semi-stable if and only if $\la \chi, \lambda\ra\ge 0$
holds for any 1-PS $\lambda$ with
$\la \lambda , \mathcal{D}_{x} \ra \ge 0$.
Therefore the set of such characters $\chi$ is the double dual cone of the cone $\mathcal{D}_{x}$.
Since $\mathcal{D}_{x}$ is rational polyhedral, the double dual coincides with
itself by \cite[(1) on page 9]{Fulton_ITV}.
Stable case can be checked similarly.
\end{proof}

Consider the effective cone $C^{T}(V)=\bigcup_{x\in V}\mathcal{D}_x$ of characters with non-empty semi-stable loci. 
We define the following notions in accordance with \cite{Dolgachev-Hu}.
\begin{definition}
A \emph{wall} defined by $x\in V$ is the set $\partial\mathcal{D}_x$.
A \emph{GIT chamber} is a connected component of the set
$C\setminus\bigcup_{x\in V}\partial\mathcal{D}_x$.
Two characters $\chi,\chi'$ are said to be \emph{wall equivalent} if the strictly semi-stable loci
$V^{\strictlysemistable}(\chi)$ and $V^{\strictlysemistable}(\chi')$ coincide. 
A connected component
of a wall equivalence class, if it is not a chamber, is called
a \emph{(GIT) cell}.
The characters $\chi,\chi'$ are said to be \emph{GIT equivalent}
if
$
 V ^{ \semistable } ( \chi )
 =
 V ^{ \semistable } ( \chi' )
$.
\end{definition}

Via similar arguments as in \cite[Theorem 3.3.2]{Dolgachev-Hu} and \cite[Lemma 3.3.10]{Dolgachev-Hu},
we can check the following statements.
\begin{lemma}\label{lm:VGIT lemmas}
\begin{enumerate}
\item A GIT chamber is a GIT equivalence class.
\item Any cell is contained in a GIT equivalence class.\label{it:cell to GIT equivalence}
\item For any GIT chamber $\cC$ we have
\begin{equation*}
\cC=\bigcap_{x\in V^{\semistable}(\cC)}\mathcal{D}_{x}^{\circ}.
\end{equation*}

\end{enumerate}
\end{lemma}


\subsection{Strong Mori equivalence is equivalent to GIT equivalence}
\label{sc:Strong Mori equivalence is equivalent to GIT equivalence}

Let $X$ be a Mori dream space, and fix a Cox ring $R=R(X, \Gamma)$.
We prove that the following three objects coincide:
\begin{itemize}
\item relative interiors of cones of $\Fan{(X)}$,
\item strong Mori equivalence classes, and
\item GIT equivalence classes.
\end{itemize}
This was proved in \cite{MR1786494} for the interiors of Mori chambers, but we need it
for all strong Mori equivalence classes in the proofs of \pref{th:fans}.
This comes from the fact that the interior of a chamber of the target space is not contained in a chamber but
a cell of the source space in general (see \pref{eg:bl-up_of_P^3}; here
the ample cone of the target space coincides with the relative interior
of a low-dimensional cell).

We first recall the VGIT of Cox rings (see \cite{MR1786494} for detail).
Consider the affine variety	 $V=\Spec{R(X, \Gamma)}$ and the torus
$T:=\Hom_{\gp}(\Gamma, \bG_m)$.
The $\Gamma$-grading of $R(X, \Gamma)$ yields an action
of $T$ on $R(X, \Gamma)$ (hence on $V_X$) in the following way;
an element $g\in T$ acts on $f\in H^0(X,\cO_X{(D)})\subset R(X, \Gamma)$
by $f\mapsto g(D)f$.

Note that there is the canonical isomorphism
\begin{equation}\label{eq:double dual}
 \ev _{ \bullet } \colon \Gamma \simto \chi { ( \chi { ( \Gamma ) } ) }
 = \chi { ( T ) },
\end{equation}
which maps $D$ to the character
$
 \ev _{D} \colon g \mapsto g ( D )
$.
After taking the tensor product $\otimes_{\bZ}\bR$,
this gives the inverse of the map
$
 \psi \colon \chi { ( T ) } _{ \bR } \simto \Gamma _{ \bR } \simto \Pic { ( X ) } _{ \bR }
$
of \cite[Theorem 2.3]{MR1786494}.

Using the isomorphism \pref{eq:double dual}, we can show that
for each $D\in\Gamma$ the ring of $\ev_D$ semi-invariants of $R(X, \Gamma)$
coincides with $R(X, \cO_X(D))$ as graded algebras.
Hence we see $Q_{\ev_D}=\Proj R(X, \cO_X(D))$, and
in particular $Q_{\ev_A}=X$ for any ample $A\in\Gamma$.
The universal property of categorical quotients provides us the rational map
\begin{equation*}
X=Q_{\ev_A}\dasharrow Q_{\ev_D},
\end{equation*}
and this coincides with the rational map
$\varphi_{\cO_X(D)}\colon X\dasharrow\Proj{R(X, \cO_X(D))}$.
Summing up, we obtain the following commutative diagram.

\[
\xymatrix{
V^{\semistable}(\ev_A)
\ar@{->>}[dd]_{/T}
&
V^{\semistable}(\ev_A) \bigcap V^{\semistable}(\ev_D)
\ar@{_{(}->}[l] \ar@{->>}[dd]_{/T} \ar@{^{(}->}[r]
&
V^{\semistable}(\ev_D) \ar@{->>}[dd]_{//T}
\\
&
&
\\
V^{\semistable}(\ev_A)/T \ar@{=}[d]
&
V^{\semistable}(\ev_A)\bigcap V^{\semistable}(\ev_D)/T
\ar@{_{(}->}[l]
\ar[r]
&
V^{\semistable}(\ev_D)//T \ar@{=}[d]
\\ 
X \ar@{-->}[rr]_{\varphi_D}
&
&
\Proj R(X, \cO_X(D))
}
\]
The three arrows
$
\hookrightarrow
$
are open immersions and $/T$ (resp. $//T$)
indicates the geometric (resp. categorical) quotient by the torus $T$.

As a consequence of these observations, we can prove the coincidence of the three notions.
\begin{proposition}\label{pr:strong Mori=GIT}
For two line bundles $L, M$ on a Mori dream space $X$, the followings are equivalent.
\begin{enumerate}
\item $L$ and $M$ are GIT equivalent; namely, $U_{\ev_{L}}=U_{\ev_{M}}$.
\item $L$ and $M$ are strongly equivalent.
\item There exists a cone $\sigma\in\Fan{(X)}$ such that
$L$ and $M$ are contained in the relative interior $\sigma^{\relint}$ of $\sigma$.
\end{enumerate}
\end{proposition}

\begin{proof}
The equivalence of the last two conditions is \pref{pr:fan vs strong Mori equivalence}.
For the equivalence of the first two, the arguments in the proof of \cite[Theorem 2.3]{MR1786494} works
essentially without change.
In the proof they only proved that the relative interiors of the Mori chambers are
identified (via $\psi$) with the GIT chambers, but the arguments can be applied
more generally to arbitrary strong Mori equivalence classes. 
We give a sketch of the proof (see the proof of \cite[Theorem 2.3]{MR1786494} for detail).

Fix a character $\chi$ which corresponds to an ample line bundle on $X$.
For an arbitrary character $y\in C^{T}(V)\cap\chi{(T)}$, let
$\psi{(y)}=P+N$ be a Zariski decomposition of the corresponding $\bQ$-line bundle.
Then we can show the equality
\begin{equation}\label{eq:2.3.2}
U_{\chi}\setminus U_{y}=q_{\chi}^{-1}(\Supp{(N)})
\end{equation}
in codimension one.
This follows from the equality
$
 H^0(X,\cO_X(m\psi{(y)}))=H^0(U_{\chi},L_{y}^{\otimes m})^{T},
$
which is (2.3.2) in the proof of \cite[Theorem 2.3]{MR1786494}. From this
we can immediately conclude that GIT equivalence implies the strong Mori equivalence.
Conversely if we assume the strong Mori equivalence of $\psi{(y)}$ and $\psi{(z)}$ for two characters $y$ and $z$,
then we see that $Q_y=Q_z$ and that $U_y$ and $U_z$ coincide in codimension one by
\pref{eq:2.3.2}. The rest of the proof
is exactly the same as that of \cite[Theorem 2.3]{MR1786494}.
\end{proof}

\begin{corollary}\label{cr:cell=GIT equivalence}
A GIT equivalence class is a GIT cell
and vice versa in the VGIT of Cox rings.
\end{corollary}
\begin{proof}
Since we have \pref{lm:VGIT lemmas} (\pref{it:cell to GIT equivalence}),
it is enough to show that any GIT equivalence class is contained in a GIT cell.
Take $\sigma\in\Fan{(X)}$. If $\sigma^{\relint}$ is not contained in a cell, the stable loci are not constant on it: i.e. there exists a point $x\in V$ such that
$\mathcal{D}_x^{\circ}\cap\sigma^{\relint}\not=\emptyset$ but $\sigma^{\relint}\not\subset\mathcal{D}_x^{\circ}$.
Since $\mathcal{D}_x$ and $\sigma$ are rational polyhedral cones, this means
$\sigma^{\relint}\not\subset\mathcal{D}_x$, contradicting the fact that $\sigma^{\relint}$ is a GIT equivalence class.
\end{proof}


\section{Comparison of the fans via GIT}\label{sc:Comparison of the fans via GIT}

In this section we prove \pref{th:fans} via the GIT interpretation of the relative interiors of the cones
(see \pref{pr:strong Mori=GIT}).
Let $f\colon X\to Y$ be a surjective morphism between Mori dream spaces. Fix a group $\Gamma_Y\subset\Divi{(Y)}$
of Cartier divisors such that the natural map $\Gamma_{\bQ}\to\Pic{(Y)}_{\bQ}$ is an isomorphism.
Also, fix a group $\Gamma_{X}\subset\Divi{(X)}$ of Cartier divisors with a similar property and which
contains $f^{*}\Gamma_Y$.
For such a pair $(\Gamma_Y, \Gamma_{X})$, the natural ring homomorphism $f^{*}\colon R(Y, \Gamma_Y)
\to R(X, \Gamma_{X})$ induces a morphism of affine varieties
\begin{equation*}
V_f\colon V_X=\Spec R(X, \Gamma_{X})\to V_Y=\Spec R(Y, \Gamma_Y).
\end{equation*}
Set $T_X=\Hom(\Gamma_{X},\bG_m)$ and $T_Y=\Hom(\Gamma,\bG_m)$.
The inclusion $f^{*}\Gamma_Y\subset\Gamma_{X}$
yields a surjective homomorphism of tori $T_f\colon T_X\to T_Y$.
Note that the morphism $V_f$ is equivariant
with respect to the actions of $T_X$ on $V_X$ and $T_Y$ on $V_Y$, and the homomorphism $T_f$.

The following is the main ingredient of the proof of \pref{th:fans}:
\begin{proposition}\label{pr:ss loci and pull back}
Let $f\colon X\to Y$ be as in \pref{th:fans}.
If we choose an appropriate pair
$(\Gamma_Y, \Gamma_{X})$ as above, for any divisor $D\in\Gamma_Y$ we have the equality
$
 V_f^{-1}(V_Y^{\semistable}(\ev_D))=V_X^{\semistable}(\ev_{f^{*}D})
$.
\end{proposition}

\begin{remark}
In general, the equality
$
 V_f^{-1}(V_Y^{\strictlysemistable}(\ev_D))=V_X^{\strictlysemistable}(\ev_{f^{*}D})
$
does not hold ($\strictlysemistable$ stands for `strictly semi-stable').
In fact, look at \pref{eg:bl-up_of_P^3} and take an ample divisor
$D$
on $Y$. Then we know that
$V_Y^{\strictlysemistable}(\ev_D)=\emptyset$.
On the other hand, since $f^{*}D$ is on a wall,
we know that
$V_X^{\strictlysemistable}(\ev_{f^{*}D})\not
=
\emptyset$.

The conclusion of \pref{pr:ss loci and pull back} does not hold for an arbitrary equivariant morphism
between affine varieties. For example, consider the morphism
\begin{equation*}
\varphi \colon \bA^{2}\to\bA^{1};(x_1,x_2)\mapsto x_1,
\end{equation*}
together with actions of $\bG_m$ on both sides with weights one. Let $\chi$ be the character of weight one.
Then
\begin{equation*}
(\bA^{2})^{\semistable}(\chi)=\bA^{2}\setminus\{0\}\supsetneq(\bA^{1}\setminus\{0\})\times\bA^{1}
=\varphi^{-1}((\bA^{1})^{\semistable}(\chi)).
\end{equation*}
\end{remark}

The following is the GIT counterpart of \pref{th:comparing strong Mori equivalence}
\begin{corollary}\label{cr:comparing GIT equivalence}
With the same assumptions as above, let $D,E$ be Cartier divisors on $Y$.
Then
$
 V _Y ^{ \semistable } ( \ev _D ) = V _Y ^{ \semistable } ( \ev _E )
$
if and only if
$
 V _X ^{ \semistable } ( \ev _{ f ^* D } )
 =
 V _X ^{ \semistable } ( \ev _{ f ^* E } )
$.
\end{corollary}
\begin{proof}
This follows from \pref{pr:ss loci and pull back} and the surjectivity of $V_f$.
\end{proof}

\noindent
\pref{th:fans} immediately follows from
\pref{cr:comparing GIT equivalence}, in view of \pref{pr:strong Mori=GIT}
(see also the proof of \pref{th:fans} in \pref{sc:Comparison of the fans without GIT}).

In the rest of this section we give a proof of \pref{pr:ss loci and pull back}.
The following lemma is the key ingredient.
\begin{lemma}\label{lm:stable loci under finite morphism}
Let $G$ be a reductive group acting on affine varieties $Z$ and $W$.
Assume that $\pi\colon Z\to W$ is an equivariant finite morphism.
Then for any $G$-linearization $\mathcal{L}$ on $W$ we have
\begin{equation*}
\pi^{-1}(W^{\semistable}(\mathcal{L}))=Z^{\semistable}(\pi^{*}\mathcal{L}).
\end{equation*}
\end{lemma}
\begin{proof}
See the proof of \cite[Theorem 1.19]{Mumford-Fogarty-Kirwan}, and
\cite[Appendix to Chapter 1, Section C]{Mumford-Fogarty-Kirwan} for
positive characteristic cases.
\end{proof}

\begin{proof}[Proof for \pref{pr:ss loci and pull back}]
Take the Stein factorization
$
 X \xxto[]{g}\Ytilde\xxto[]{h}Y.
$
Set $V_{\Ytilde}=\Spec R(\Ytilde, h^*\Gamma_Y)$.
Since $h^*\colon \Gamma_Y\to h^*\Gamma_Y$ is an isomorphism,
$T_Y$ acts on $V_{\Ytilde}$ and the natural map
$V_h\colon V_{\Ytilde}\to V_Y$ is equivariant.
Since $V_h$ is finite by \pref{th:finite},
\pref{lm:stable loci under finite morphism} implies the equality
$
 V_{\Ytilde}^{\semistable}(\ev_{h^*D})=V_h^{-1}(V_Y^{\semistable}(\ev_{D})).
$
Therefore it is enough to show
$
 V_X^{\semistable}(\ev_{f^*D})=V_g^{-1}(V_{\Ytilde}^{\semistable}(\ev_{h^*D})).
$

Since $V_g$ is an affine morphism,
the inclusion $\supseteq$ directly follows from the definition of semi-stability.
Conversely take a point $x\in V_X^{\semistable}(ev_{f^*D})$. By the definition of semi-stability,
there exists an element $s\in R(X, \Gamma_X)$ such that
$s(x)\not=0$ and $g\cdot s=g(D)^ms$ holds for any $g\in T_X$.
It is easy to check that this is equivalent to the condition $s\in R(X, f^*D)$.
Since $g$ is an algebraic fiber space, $s$ is the pull-back of an section of $R(\Ytilde, h^*D)$.
This concludes the proof.
\end{proof}


\section{Examples}\label{sc:Examples}

\begin{example}\label{eg:bl-up_of_P^3}
We borrow from \cite[Example 5.5]{MR2825271}.
Let $X$ be the blow-up of $\bP^3$ in two distinct points, say $p_1$ and $p_2$.
Since $X$ is toric, it is a Mori dream space. Let 
$E_i$ be the exceptional divisor over $p_i $ $( i = 1, 2 )$, and
$\ell$ the line passing through the points $p_1,p_2$.
Let $E_3$ be the class of the strict transform of a plane containing $\ell$.
We can show that $X$ has a flopping contraction which contracts
the strict transform of the line $\ell$. Let $X'$ be the flop. Using toric descriptions,
we see that this is an Atiyah flop.

The effective cone of $X$ is spanned by the divisors $E_i$, and
the movable cone is the union of the semi-ample cones of $X$ and $X'$.
$\Nef{(X)}$ is spanned by three divisors $H,H-E_1$, and $H-E_2$, where
$H$ is the total transform of a plane of $\bP^3$. $\Nef{(X')}$ is
spanned by $H-E_1,H-E_2$, and $E_3$.

A slice of $\Eff{(X)}$, together with its fan structure is described
in the following figure:

\vspace{5mm}

$
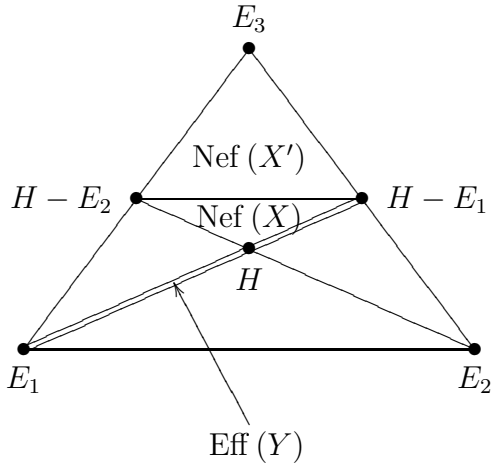
\begin{xy}
(0,0)="A"*{\bullet},
"A"+<6cm,0cm>="B"*{\bullet},
"A"+<3cm,4cm>="C"*{\bullet},
"A"+<1.5cm,2cm>="D"*{\bullet},
"D"+<3cm,0cm>="E"*{\bullet},
"A"+<3cm,1.33333cm>="F"*{\bullet},
(0,-4)="E_1"*{E_1},
"B"+(0,-4)="E_2"*{E_2},
"C"+(0,4)="E_3"*{E_3},
"D"+(-10,0)="H-E_2"*{H-E_2},
"E"+(10,0)="H-E_1"*{H-E_1},
"F"+(0,-4)="H"*{H},
"F"+(0,4)="SA(X)"*{\Nef{(X)}},
"F"+(0,12)="SA(X')"*{\Nef{(X')}},
"A"+<2cm,0.888889cm>="T",
"A"+<3cm,-1cm>="S",
"S"+(1,-3)*{\Eff{(Y)}}
\ar"S";"T"
\ar@{-}"A";"B"
\ar@{-}"A";"C"
\ar@{-}"B";"C"
\ar@{=}"A";"E"
\ar@{-}"B";"D"
\ar@{-}"D";"E"
\end{xy}$
\vspace{5mm}

Let $Y$ be the blow-up of $\bP^3$ in $p_1$. Then the effective cone of $Y$, together with
its fan structure sits in $\Eff{(X)}$ as indicated in the figure above (the slice of $\Eff{(Y)}$ is
denoted by the double line).
As indicated in the figure above, $\Eff{(Y)}$ is mapped onto the cone spanned by
$E_1$ and $H-E_1$. The cone spanned by $H$ and $H-E_1$ is the semi-ample cone of $Y$, and
the one spanned by $H$ and $E_1$ corresponds to the Mori chamber of $Y$ whose interior
correspond to the line bundles whose positive part defines the birational contraction to $\bP^3$
and the support of whose negative part is the exceptional divisor of the contraction.

Now take a coordinate on $\bP^3$ and assume by a linear coordinate change that $p_1=(0:0:0:1)$ and $p_2=(0:0:1:0)$. Consider the
action of $\bZ_2$ on $\bP^3$ defined by $(x:y:z:w)\mapsto (x:y:w:z)$. This action lifts to
$X$, and let $X\to Z$ be the quotient
morphism. The effective cone of $Z$ together with its Mori chamber decomposition
sits in that of $X$ as follows ($\Eff{(Z)}$ is
denoted by the double line):

\vspace{5mm}
$
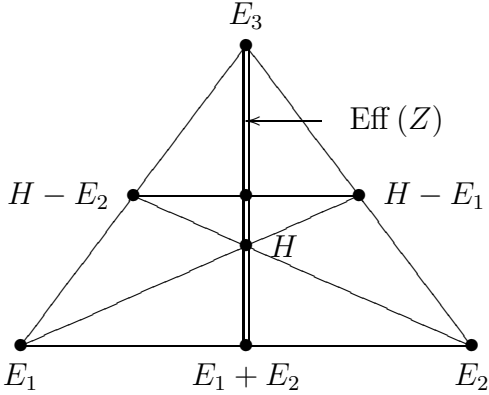
\begin{xy}
(0,0)="A"*{\bullet},
"A"+<6cm,0cm>="B"*{\bullet},
"A"+<3cm,4cm>="C"*{\bullet},
"A"+<1.5cm,2cm>="D"*{\bullet},
"D"+<3cm,0cm>="E"*{\bullet},
"A"+<3cm,1.33333cm>="F"*{\bullet},
"A"+<3cm,0cm>="G"*{\bullet},
"A"+<3cm,2cm>="I"*{\bullet},
(0,-4)="E_1"*{E_1},
"B"+(0,-4)="E_2"*{E_2},
"C"+(0,4)="E_3"*{E_3},
"D"+(-10,0)="H-E_2"*{H-E_2},
"E"+(10,0)="H-E_1"*{H-E_1},
"F"+(5,0)="H"*{H},
"G"+(0,-4)="E_1+E_2"*{E_1+E_2},
"C"+<0cm,-1cm>="T",
"T"+<1cm,0cm>="S",
"S"+(10,0)="Eff{(Z)}"*{\Eff{(Z)}},
\ar"S";"T"
\ar@{-}"A";"B"
\ar@{-}"A";"C"
\ar@{-}"B";"C"
\ar@{-}"A";"E"
\ar@{-}"B";"D"
\ar@{-}"D";"E"
\ar@{=}"G";"C"
\end{xy}$
\vspace{5mm}

As indicated in the diagram above, we can see that $Z$ has two Mori chambers other than
the semi-ample cone (recall that $\Nef{(Z)}$ coincides with the restriction of $\Nef{(X)}$ to
$\Pic{(Z)}_{\bR}$).

Let $Z'$ be the quotient of $X'$ by the involution induced from that on $X$.
Note that the Mori chamber of $Z$ obtained by restricting $\Nef{(X')}$ is the
semi-ample cone of $Z'$. The morphism defined by the ray separating $\Nef{(Z)}$ and
$\Nef{(Z')}$ is the flipping contraction of $Z$ which contracts the image of $\ell$ under the quotient
morphism $X\to Z$, and $Z'$ is the flip. 

This example shows that a Mori chamber of the target space $Z$ is not necessarily a face of
a Mori chamber of the source $X$.
\end{example}


\begin{example}\label{eg:globally F-regular is not preserved}
This example is well known to experts, but we give a detailed explanation for the sake of completeness.
The author learned this example from Tadakazu Sawada.

Suppose that $\bfk=\overline{\bfk}$ and $\chara{\bfk}=p>0$.
Let $\bA^{2}_{x,y}\subset\bP^2$ be a standard embedding of affine $2$-plane
with coordinates $x$ and $y$. Take $f=f(x,y)\in \bfk[x,y]$.
Consider the rational vector field defined by
\begin{equation*}
\delta=\frac{\partial f}{\partial y}\frac{\partial}{\partial x}-\frac{\partial f}{\partial x}\frac{\partial}{\partial y}
\end{equation*}
and the $p$-closed foliation on
$\bP^2$ it generates.
By the Jacobson correspondence described in \pref{pr:Jacobson correspondence}
we obtain the quotient morphism
$\pi\colon \bP^2\to Y$,
which is a purely inseparable finite morphism of degree $p$
to a normal projective variety $Y$.
The coordinate ring of the open subset
$\pi{(\bA^2)} \subset Y$
is isomorphic to
\begin{align*}
 \bfk[x,y] ^{\delta}
 =
 \bfk[x^p, y^p, f(x,y) ]
 \simeq
 \bfk[X,Y,Z]/(Z^p-f(X,Y)).
\end{align*}

Set $f(x,y)=x^py+xy^p$. By Fedder's criterion for $F$-purity (\cite[Proposition 1.7]{MR701505}),
we can check that
the singularity $(0\in \bfk[x,y]^{\delta})\simeq \bfk[[X,Y,Z]]/(Z^p-X^pY-XY^p)$ is not $F$-pure.
Therefore $Y$ is not globally $F$-regular, even though $\bP^2$ is globally $F$-regular.

On the other hand we can show that $Y$ is a Mori dream space.
Firstly the Picard number of $Y$ is one since $\pi^{*}\colon \Pic{(Y)}_{\bR}\to \Pic{(\bP^{2})}_{\bR}$ is
injective. Together with the claim below, by \pref{lm:pic}, we see that $Y$ is a Mori dream space
of Picard number one.
\begin{claim}
Let $f\colon X\to Y$ be a purely inseparable finite morphism between normal varieties. If $X$ is $\bQ$-factorial, so is
$Y$.
\end{claim}

\begin{proof}
Let $(U,f_U)_U$ be a Cartier divisor on $X$, where $X=\bigcup U$ is an
open cover of $X$ and $f_U\in \bfk(X)$. Then we can check that the pushforward
$f_*(U,f_U)_U$ as a Weil divisor coincides with the Cartier divisor $(U,N(f_U))_{U}$,
where
$
 N = N _{ X / Y } \colon \bfk(X) \to \bfk ( Y )
$
is the norm function.

Now let $D$ be a Weil divisor on $Y$. By assumption, there exists a positive integer $m$ such that
$mf^{*}D=f^{*}mD$ is Cartier. Since
$
 f _{*} f ^{*} m D 
$
is Cartier and it coincides with
$
 m \deg{ ( f ) } D
$,
we get the conclusion.
\end{proof}
\end{example}

\begin{remark}
Recall that the notions of
$F$-purity and global $F$-regularity are conjecturally
related, via reductions to positive characteristics,
to log terminality and varieties of Fano type in characteristic zero
(see \cite{MR2628797}).
Since images of varieties of Fano type in characteristic zero are
always of Fano type by \cite[Corollary 5.2]{MR2944479},
this example might seem to disprove the conjecture.
This is not the case, since the precise conjecture asserts that
being of Fano type should be equivalent to being of globally
$F$-regular \emph{for almost all primes}:
in the example above, we can check that
the lift of the quotient
$Y$
to characteristic zero
reduces to a Fano variety in every other characteristic than
$p$. 
\end{remark}

\begin{example}\label{eg:failure_of_being_MDS_under_operations}
Here we collect several known examples to show that
the property of being a Mori dream space is not preserved by
standard operations on varieties.
We also mention the behavior of the fan structure.

\begin{itemize}
\item
Blowup of
$
\bP^2
$
in a general set of eight points is a del Pezzo surface, and hence
is a Mori dream surface.
On the other hand, blowing-up in general nine points is not a Mori dream space.
The property of being a Mori dream space is not a birational invariant.

\item Even worse, the crepant resolution of a Mori dream space is not necessarily a Mori dream space.
Consider K3 surfaces $S$ of Picard number $20$
whose canonical models have Picard number one.

The canonical model is a Mori dream surface for trivial reasons.
We can show that any K3 surface of Picard number $20$ is not a Mori dream space,
since its automorphism group is always discrete and infinite (see \cite{MR0441982})
and by \cite[Theorem 2.11]{MR2660680}.
See \cite[Examples 1,2]{MR1420924} for examples of such surfaces $S$.

\item Blowing-up of $\bP^2$ in a set of points lying on a line is always a Mori dream surface
by \cite[Example 3.3]{MR2058459}. On the other hand, if the number of points is at least nine
and we let the points move out of the line, the blowing-up fails to be a Mori dream surface in general.
Therefore the property of being a Mori dream space is not an open condition.

\item Projective space bundles over Mori dream spaces are not necessarily Mori dream spaces (see \cite{MR2968631}).

\item It is known that for a family of terminal $\bQ$-factorial Fano varieties,
the effective cones and the movable cones of the fibers are locally constant
\cite[Corollary 5.1, Theorem 5.6]{MR2931867}. Moreover,
the Cox rings of the fibers form a flat family \cite[Corollary 4.5]{MR2931867}.
On the other hand, the Nef cone might change in such a family in general
\cite[Theorem 1.1]{MR2877439}.
\end{itemize}

\end{example}

\section{Amplifications}\label{sc:Amplifications}

In this section, we extend our results to
\begin{itemize}
\item the non-$\bQ$-factorial case, and
\item Mori dream regions.
\end{itemize}

\subsection{Mori dream space without the $\bQ$-factoriality condition}\label{sc:Non-Q-factorial Mori dream space}
In \cite[Section 2]{MR2660680}, the notion of Mori dream space was extended to normal projective varieties which are not
necessarily $\bQ$-factorial. Throughout this subsection we call them Mori dream spaces, and
use the term $\bQ$-factorial Mori dream space to indicate the ordinary one of \pref{df:Mori dream space}.
We can show that our main results are also valid for this Mori dream space.
For simplicity we assume that the base field is algebraically closed.

\begin{definition}\label{df:nQfMori dream space}
Let $X$ be a normal projective variety whose divisor class group
$\Cl{(X)}$ is finitely generated.
\begin{enumerate}
\item A Cox ring of $X$ is a multi-section ring $R(X, \Gamma)$, where
$\Gamma\subset \WDivi{(X)}$ is a group of Weil divisors such that the natural map
$\Gamma_{\bQ}\to\Cl{(X)}_{\bQ}$ is an isomorphism.
\item $X$ is said to be a Mori dream space if a Cox ring of $X$ is of finite type over the base field.
\end{enumerate}
\end{definition}
Note that Mori dream spaces and Cox rings in the above sense coincide with those of
$\bQ$-factorial Mori dream spaces when the variety is $\bQ$-factorial.
The finite generation of Cox rings is independent of the choice of $\Gamma$.

In \cite[Theorem 2.3]{MR2660680}, they gave a characterization of
Mori dream spaces via the properties of line bundles,
which is similar to the original definition.
In \cite[Theorem 2.3]{MR2660680} and its proof they also proved the existence of small $\bQ$-factorizations
for such varieties, assuming that the base field is of characteristic zero.
They need this extra assumption in the proof of \cite[Lemma 2.4]{MR2660680},
to use resolutions of singularities.
We prove that it is unnecessary, so that their results work in arbitrary characteristics.

\begin{lemma}[{$=$\cite[Lemma 2.4]{MR2660680}}]
Let $X$ be a normal projective variety with finitely generated class group $\Cl{(X)}$.
Then the cone $\Mov{(X)}$ generated by divisors without fixed parts
has full dimension in $\Cl{(X)}_{\bQ}$.
\end{lemma} 

\begin{proof}
Take prime divisors $D_1,\dots,D_r$ on $X$ which generates $\Cl{(X)}$.
Let $f\colon X'\to X$ be the
composition of normalizations and successive blow-ups along scheme-theoretic
inverse images of $D_i$.
Let $D_i'$ be the scheme theoretic inverse image of $D_i$ on $X'$,
which are Cartier by the construction of the blow-ups and satisfy $f_*D'_i=D_i$.
Let $E$ be a very ample Cartier divisor on $X'$ such that $E+D'_i$ are also very ample.
Since $f_*E$ and $f_*(E+D'_i)=f_*E+D_i$ are all movable on $X$, this concludes the proof.
\end{proof}

We go back to our results. We first show \pref{th:main} for non-$\bQ$-factorial varieties.
\begin{theorem}\label{th:main'}
Let $X$ be a Mori dream space, and $X\to Y$ be a surjective morphism to another
normal projective variety. Then $Y$ also is a Mori dream space.
\end{theorem}
\begin{proof}
We only point out which part of the proof of \pref{th:main} should be modified.
Let us first replace $X$ with its small $\bQ$-factorization. Namely, take
a small birational morphism $\tilde{X}\to X$ from a $\bQ$-factorial Mori dream space
$\tilde{X}$. By replacing $X$ with $\tilde{X}$, we can assume that $X$ is $\bQ$-factorial.

In order to prove the finite generation of a Cox ring of $Y$, we take the Stein factorization of $X\to Y$
as before. Nothing has to be changed for finite morphisms.

Suppose that $f$ is an algebraic fiber space. Here we need some non-trivial
modifications, and we quote some arguments from the proof of \cite[Theorem 1.2]{MR2811268}.
Denote by
$
 U \subset Y
$
the non-singular locus of
$Y$ and set
$
 V = f^{-1} ( U )
$.
Let
$
 E_1, \dots, E_m
$
be the codimension one irreducible components of
$
 X \setminus V
$
and set
$
 E = \sum _i E _i
$.
Fix a global section
$
 f \in H^0 ( X, \cO _X ( E ) )
$.
Take any group of Weil divisors $\Gamma$ on $Y$ and
$
 \Gammabar
$
on
$
 X
$
such that
$
 ( f | _V ) ^* ( \Gamma | _U )
 =
 \Gammabar | _V
$.
Set
$
 \Gammatilde = \Gammabar \oplus\bZ E
$.
Then we can show the isomorphism
$
 R ( X, \Gammatilde ) _f \simeq R ( V, \Gammatilde | _V ).
$
On the other hand since the restriction of $\cO_X(E)$ to $V$ is trivial, we see that
$
 R ( V, \Gammatilde | _V )
$
is isomorphic to
$
 R \lb V, ( f | _V ) ^* ( \Gamma | _U ) \rb [ t ^{ \pm 1 } ]
 \simeq
 R ( Y, \Gamma ) [ t ^{ \pm 1 } ]
$.
Now we can apply \pref{lm:afs} to $X$ and $\Gammatilde$
to obtain the finite generation of
$
 R ( X, \Gammatilde )
$.
From the isomorphism we have just seen,
the finite generation of $R(Y,\Gamma)$ follows. Thus we conclude the proof.
\end{proof}

For a Mori dream space $X$, we can define the notion of Mori equivalence, Mori chambers and so on
as those of its small $\bQ$-factorial modifications. The choice of the modifications do not matter.
In particular we can define the fan of $X$, and for this we can show the result generalizing
\pref{th:fans}.

\begin{theorem}\label{th:fans'}
Let $X\to Y$ be a surjective morphism between Mori dream spaces. Then
\begin{equation*}
\Fan{(Y)}=\Fan{(X)}|_{\Eff{(Y)}}.
\end{equation*}
\end{theorem}
\begin{proof}
By taking suitable small $\bQ$-factorizations of $X$ and $Y$, the morphism lifts to the one between
$\bQ$-factorial Mori dream spaces. Thus we can reduce the problem to our original \pref{th:fans}.
\end{proof}

%
%

\subsection{Mori dream region}\label{sc:Mori dream region}

Let $X$ be a normal $\bQ$-factorial projective variety.
There is a notion called Mori dream region defined in \cite[Definition 2.12]{MR1786494},
which generalizes Mori dream spaces. In this subsection we check that \pref{th:main} can be extended to
Mori dream regions.

First we recall the definition of Mori dream regions from \cite[Definition 2.12]{MR1786494}:
\begin{definition}\label{df:Mori dream region}
Let $X$ be a normal $\bQ$-factorial projective variety. A cone $\cC\subset\Pic{(X)}_{\bR}$
spanned by finitely many line bundles is called a Mori dream region if
the multi-section ring
\begin{equation*}
R(X, \cC) = \bigoplus _{ D \in \cC \cap \Pic { (X)} ^{\free} } H^{0} ( X, \cO _{X} ( D ) )
\end{equation*}
is of finite type over the base field.
\end{definition}

\noindent
If the natural morphism $\Pic{(X)}_{\bQ}\to\Num{(X)}_{\bQ}$ is an isomorphism,
by definition, $\Eff{(X)}$ is a Mori dream region if and only if $X$ is a Mori dream space.

As in the case of Mori dream space, Mori dream region can be characterized via the existence of a decomposition into
finitely many rational polyhedral subcones such that on each of them Zariski decomposition is
$\bQ$-linear:

\begin{proposition}\label{pr:Mori dream region and ZD}
Let $X$ be a normal $\bQ$-factorial projective variety and
$
 \cC \subset \Pic { ( X ) } _{ \bR }
$
a cone spanned by finitely many line bundles.
It is a Mori dream region if and only if the following conditions are satisfied:
\begin{itemize}
\item
$
 \cC \cap \Eff { ( X ) }
$
is spanned by finitely many line bundles.
\item
The section ring of any line bundle in
$
 \cC \cap \Eff { ( X ) }
$
is of finite type over the base field.
In particular, $\bQ$-effective line bundles admit Zariski decompositions:
i.e., for such a line bundle $L$ there exists a decomposition $L=P+N$ such that
$P$ is movable and $N$ is effective, and for all sufficiently divisible positive integer $m$
all the global sections of $mL$ come from those of $mP$ as in \pref{eq:ZD} of \pref{pr:ZD}.

\item There exists a decomposition of
$
 \cC \subset \Pic { ( X ) } _{ \bR }
$
into finitely many rational polyhedral subcones such that Zariski decomposition is $\bQ$-linear
on each of them.
\end{itemize}
\end{proposition}
\begin{proof}
`if' part is exactly the same as the proof of \pref{lm:afs}.
For the `only if' part, the second condition follows from \cite[`If' part of Lemma 1.6]{MR1786494}. For the third condition,
see \cite[Theorem 3 (3)]{MR3048609}.
\end{proof}

\begin{remark}
In \cite[Theorem 2.13]{MR1786494}, they claim that we can find a decomposition of $\cC$ into
chambers $\cC_i$ and for each of them
there exists a contracting birational map
$g_i\colon X\dasharrow Y_i$ such that
\begin{equation*}
 \cC _i = \cC \cap \lb g _i ^{*} \Nef { ( Y _i ) } *\ex { ( g _i ) } \rb.
\end{equation*}

The author believes that it is not so easy to prove, since it means that
the existence of the canonical model would imply the existence of a minimal model.
This is why he replaced \cite[Theorem 2.13]{MR1786494} with \pref{pr:Mori dream region and ZD}.

\end{remark}

\begin{corollary}\label{cr:main''}
Let $f\colon X\to Y$ be a surjective morphism between normal $\bQ$-factorial projective varieties.
Let
$
 \cC \subset \Pic { ( X ) } _{ \bR }
$
be a finitely generated rational polyhedral cone which is a Mori dream region. Then
$\cC|_{\Eff{(Y)}}$ also is a Mori dream region.
\end{corollary}
The proof is essentially the same as the one for \pref{th:main}.

\begin{remark}
\pref{th:fans} does not hold for Mori dream regions in general.
Here we give some observations to this problem.

Consider a rational polyhedral cone which is contained in the ample cone of a normal
projective variety. Obviously, all the divisors in the cone are strongly Mori equivalent.
Taking this cone as a Mori dream region, we see that we do not have the fan structure for Mori dream regions
such that the relative interior of the cone of the fan is an equivalence class.

Moreover, we do not know if the GIT equivalence and the strong Mori equivalence coincide for
arbitrary Mori dream regions or not. The reason is as follows.
By closely looking at the proof of \pref{pr:strong Mori=GIT}, we see that
GIT equivalence implies strong Mori-equivalence for arbitrary Mori dream regions,
provided that they contain ample divisors.
For the proof of the converse, it was essential that the unstable locus for ample divisors
has codimension at least two. Even if we take an arbitrary Mori dream region which contains an
ample divisor and consider the spec of the corresponding multi-section ring, the unstable locus
of ample divisors can have divisorial component: the Nef cone of the blow-up of $\bP^2$
at a point gives such an example. The difference comes from the fact that
\cite[Lemma 2.7]{MR1786494} holds only for Cox rings.

Nevertheless \pref{cr:comparing GIT equivalence} and
\pref{th:comparing strong Mori equivalence} holds for an arbitrary Mori dream region $\cC$
on $X$ and
$
 \cC | _{ \Eff { ( Y ) } }
$.
\end{remark}

%
%

\bibliographystyle{amsalpha}
\bibliography{mainbibs}

\noindent
Shinnosuke Okawa

Department of Mathematics,
Graduate School of Science,
Osaka University,
Machikaneyama 1-1,
Toyonaka,
Osaka,
560-0043,
Japan.

{\em e-mail address}\ : \  okawa@math.sci.osaka-u.ac.jp
\ \vspace{0mm} \\

 \end{document}